\documentclass[11pt]{article}
\usepackage{amssymb,amsmath,latexsym}
\usepackage{tablists}
\usepackage{hyperref}
\usepackage[latin1]{inputenc}
\usepackage[french,english]{babel}
\usepackage{graphicx}
\newtheorem{ex}{Example}

\newtheorem{remark}{\\Remark}
\newtheorem{pro}{Proposition}
\newtheorem{lemma}{Lemma}
\newtheorem{thh}{Theorem }
\newtheorem{co}{Corollary}
\newcommand{\p}{{{\bf Proof.\,\,}}}

\newcommand{\biindice}[3]%
{

\begin{array}[t]{c}
#1\\
{\scriptstyle #2}\\
{\scriptstyle #3}
\end{array}

}

\begin{document}
\title{\Large  \textbf{Spectral Mapping Theorems of Differentiable $C_0$ Semigroups }}
\author{\normalsize Abdelaziz Tajmouati$^{1}$, Hamid Boua$^{2}$ and Mohammed Karmouni$^{3}$\\\\
\normalsize $^{(1,2)}$Sidi Mohamed Ben Abdellah University\\
\normalsize Faculty of Sciences Dhar El Mahraz \\
\normalsize Laboratory of Mathematical Analysis and Applications\\
\normalsize Fez, Morocco\\
\normalsize Email: abdelaziz.tajmouati@usmba.ac.ma\,\,\,\, hamid12boua@yahoo.com\\
\normalsize $^{(3)}$Cadi Ayyad University\\
\normalsize Multidisciplinary Faculty \\
\normalsize Safi, Morocco\\
\normalsize med89karmouni@gmail.com
}

\date{}\maketitle
\begin{abstract}
Let $(T(t))_{t\geq 0}$ be a $C_0$ semigroup on a Banach space $X$ with infinitesimal generator $A$.
 In this work, we give   conditions for which the  spectral mapping theorem $\sigma_{*}(T(t))\backslash \{0\}=\{e^{\lambda s}, \lambda\in\sigma_{*}(A)\}$
holds, where $\sigma_*$ can be equal to the essential, Browder and Kato spectrum. Also, we will be interested in the relations
between the spectrum of $A$ and the spectrum of the nth derivative $T(t)^{(n)}$ of  a differentiable $C_0$ semigroup $(T(t))_{t\geq0}$.
\end{abstract}
\section{\large Introduction and Preliminaries}

Many equations of mathematical physics can be cast in the abstract form:
\begin{equation}\label{b}
 u'(t)=Au(t), t\geq 0, u(0)=u_0
\end{equation}
 on  a Banach space $X$. Here $A$ is a given linear operator with domain $D(A)$ and the
 initial value $u_0$. The semigroups can be used to solve a large class of problems commonly known as the Cauchy problem.
 The solution of (\ref{b}) will be given by $u(t)=T(t)u_0$ for an operator semigroup $(T(t))_{t\geq 0}$ on $X$.

In order to understand the behavior of the solutions in terms of the data concerning $A$, one seeks
information about the spectrum of $T(t)$ in terms of the spectrum of $A$. Unfortunately
the spectral mapping theorem  $e^{t\sigma_*(A)}=\sigma_*(T(t))\setminus \{0\}$ often fails,  sometimes in
dramatic ways. However, the inclusion
\begin{equation}\label{a}
e^{t\sigma_*(A)}\subseteq\sigma_*(T(t))\setminus\{0\}
\end{equation}
 always holds, where $\sigma_*\in\{\sigma, \sigma_{e}, \sigma_{b}\}$  is the spectrum, essential spectrum or Browder spectrum.


A one-parameter family $(T(t))_{t\geq0}$ of bounded  operators on  a Banach space $X$ is called a $C_{0}$ semigroup of operators or a strongly continuous semigroup of operators if:
  \begin{enumerate}
    \item $T(0) = I$.
    \item $T(t + s) = T(t)T(s)$, $\forall t, s\geq 0$.
    \item $\displaystyle{\lim_{t\rightarrow 0}} T(t)x = x$, $\forall x\in X$.
  \end{enumerate}
 $(T(t))_{t\geq0}$ has a unique infinitesimal generator $A$ defined on the domain $D(A)$ by:
  $$Ax = \displaystyle{\lim_{t\rightarrow 0}} \frac{T(t)x-x}{t}, \forall x\in D(A),$$
$$D(A)=\{x\in X : \displaystyle{\lim_{t\rightarrow 0}} \frac{T(t)x-x}{t} \mbox{ exists} \}.$$
  Also,  $A$ is a closed operator, see  \cite{engel, Pazy}.\\

    We introduce the following operator acting on $X$ and depending on the parameters $\lambda \in \mathbb{C}$ and $t\geq 0$:
\begin{center}
$B_\lambda(t)x=\int_0^t e^{\lambda (t-s)}T(s)xds,  x\in X  $
\end{center}
It is well known  that $B_\lambda(t)$ is a bounded linear operator from $X$ to $D(A)$  and we have (\cite{engel, Pazy}):
\begin{description}
  \item[a. ] $(e^{\lambda t}-T(t))^nx =(\lambda -A)^nB^n_\lambda(t)x,\;\;\;   \forall x\in X,  n\in\mathbb{N};$

  \item[b. ] $(e^{\lambda t}-T(t))^nx =B^n_\lambda(t)(\lambda -A)^n x,\;\;\;   \forall x\in D(A^n),   n\in\mathbb{N}$;

  \item[c. ] $R((e^{\lambda t}-T(t))^n)\subseteq R((\lambda -A)^n)$;

  \item[d. ]   $N((\lambda -A)^n)\subseteq N(e^{\lambda t}-T(t))^n.$
\end{description}
The semigroup $(T(t))_{t\geq 0}$ is called differentiable for $t> t_0$ if for every $x\in X$,
 $t \longmapsto T(t)x$ is differentiable for $t > t_0$. $(T(t))_{t\geq 0}$ is called differentiable if it is differentiable for $t > 0$.
If $(T(t))_{t\geq 0}$ is differentiable, $B_{\lambda}(t)x$ is differentiable  and $B'_{\lambda}(t)x=T(t)x+\lambda B_{\lambda}(t)x$.\\
Let $(T(t))_{t\geq 0}$ be a $C_0$-semigroup differentiable on $X$ with infinitesimal generator $A$. We have (\cite{Pazy}):
\begin{description}
\item[a. ] $(\lambda e^{\lambda t}-AT(t))x =(\lambda -A)B'_\lambda(t)x,\;\;\;   \forall x\in X;$
\item[b. ] $(\lambda e^{\lambda t}-AT(t))x =B'_\lambda(t)(\lambda -A)x,\;\;\;   \forall x\in D(A)$;
\end{description}

 Let $\Delta =\{z\in \mathbb{C} : \alpha_1 < \mbox{ arg z } < \alpha_2 \}$, $-\frac{\pi}{2}< \alpha_1< \alpha_2<\frac{\pi}{2}$, and for $z \in  \Delta$
let $T(z)$ be a bounded linear operator. The family $(T(z))_{z \in  \Delta}$ is an analytic semigroup in $\Delta$ if
\begin{description}
\item [i. ] $z \mapsto T(z)$ is analytic in $\Delta$.
\item [ii. ] $T(0)= I$ and $\displaystyle{\lim_{z\rightarrow 0}}T(z)x = x$ for every $x \in X$.
\item [iii. ] $T(z_1 + z_2) = T(z_1)T(z_2)$ for $z_1, z_2 \in  \Delta$.
\end{description}
A semigroup $(T(t))_{t\geq 0}$ will be called analytic if it is analytic in some sector $\Delta$ containing the nonnegative real axis.

One of the classical approaches to finding information about the solution $T(t)$ is to study the spectrum of the semigroup directly.
In several applications, we have the explicit expression of the generator $A$. Thus we need information on the spectrum of the semigroup $T(t)$ in terms of that of the generator $A$.
In addition, the study of the spectral mapping theorem for the different parts of the spectrum proves to be essential.

The spectral inclusions for various reduced spectra of a differentiable $C_0$ semigroup were studied by
 A. Tajmouati et. al \cite{taj}\cite{taj2}. But,  the equality (spectral mapping theorem) is not examined. So it is natural to ask the question: can we establish equality?

In this work, we will continue in this direction,   we will prove, under some conditions, the
spectral equality for $C_0$ semigroup for essential spectrum, Browder spectrum,  upper(lower)semi-Fredholm spectrum. Also,
we give an affirmative answer to the question we just asked.
\section{$C_0$ semigroup}
Throughout, $X$ denotes a complex Banach space, let $A$ be a closed linear operator on $X$ with domain  $D(A)$, we denote by $A^*$,  $R(A)$, $N(A)$,  $\sigma_{ap}(A)$, $\sigma_{s}(A)$, $\rho(A)$ and  $\sigma(A)$,
 respectively the adjoint, the range, the null space,  the approximate point spectrum, the  surjective spectrum,  the resolvent set and   the spectrum of $A$.\\

We start by  the following results, they  will be needed in the sequel.
\begin{lemma}\label{llla}
Let $(T(t))_{t\geq 0}$ be a $C_0$ semigroup on $X$ with infinitesimal generator $A$. Suppose that there exists $t_0>0$ such that $AT(t_0)$ is bounded. Then
$(T(t)_{|F})_{t\geq 0}$ is a uniformly continuous $C_0$ semigroup on $X$, where $F=N(e^{\lambda s}-T(s))$, $\lambda\in\mathbb{C}$,  $s>0$.
\end{lemma}
\p
It is clear that $F$ is a closed subspace of $X$,  $A$-invariant  and $T(t)$-invariant, 	
furthermore  $D(A_{\mid F})=D(A)\cap F\subseteq F$. Let $x\in F$ and $x_n=n\displaystyle\int^{t+\frac{1}{n}}_t T(s)xds\in D(A)\cap F$, $x_n \underset{n \to \infty}{\longrightarrow}0$, so that $\overline{D(A_{\mid F})}=F$. Let us show that $A_{|F}$ is bounded. Let $\mu\in\rho(A)$, then $R(\mu,A_{|F})$ is bounded below. Indeed: If not, there exists $(x_n)_n\subset F$  such that $\parallel  x_n \parallel=1$ and $R(\mu,A_{|F})x_n\underset{n \to \infty}{\longrightarrow}0$. Since $AT(t_0)$ is  bounded, then
$(\mu -A_{\mid F})T(t_0)R(\mu, A_{\mid F})x_n\underset{n \to \infty}{\longrightarrow}0$. So, $T(t_0)x_n\underset{n \to \infty}{\longrightarrow}0$.\\
If $t\geq t_0$, then $T(t)x_n=T(t_0)T(t-t_0)x_n=T(t-t_0)T(t_0)x_n\rightarrow 0$, then $e^{\lambda t}x_n=T(t)x_n\rightarrow 0$, which is impossible.\\
If $0<t<t_0$, there exists  $p\in\mathbb{N}^*$ such that $pt>t_0$, then $T(pt)x_n=(T(t))^px_n=T(pt-t_0)T(t_0)x_n\rightarrow 0$, so
$e^{p\lambda t}x_n=(T(t)^p)x_n\rightarrow 0$,  then $x_n\rightarrow 0$ which contradicts the fact that $\parallel  x_n \parallel=1$.
Hence, there exists $M\in\mathbb{R}^{+}$ such that for all $x\in F$ $\parallel R(\mu, A)x\parallel\geq M\parallel x\parallel$.
Whence, for all $x\in D(A_{\mid F})$
\begin{center}
  $\parallel R(\mu, A)(A-\mu)x\parallel\geq M\parallel (A-\mu)x \parallel$.
\end{center}
Which gives
\begin{center}
    $\parallel A_{\mid F} x\parallel\leq (|\mu|+ M^{-1})\parallel x \parallel$
\end{center}
Now, let $x\in F=\overline{D(A_{\mid F})}$, there exists $(u_p)_p\subset D(A_{\mid F})$, $u_p\rightarrow x$.
For all $p,q\in\mathbb{N}$, $$\|Au_p-Au_q\|\leq (|\lambda_0|+ C^{-1})\parallel u_p-u_q \parallel$$
Hence is a Cauchy sequence on $X$, so converge. Since $A$ is a closed operator, then $x\in D(A_{\mid F})$. Hence the result.$\blacksquare$

A closed linear operator $A$ is called an upper semi-Fredholm
(resp, lower semi Fredholm) if $\dim N(A)<\infty \mbox{ and  } R(A) \mbox{  is closed }$
(resp, $codim R(A) <\infty$). $A$ is a Fredholm operator if it is a lower and upper semi-Fredholm operator,
The essential  spectrum of $T$  is defined  by :
$$\sigma_{e}(T)=\{\lambda\in \mathbb{C}:  T-\lambda I \mbox{  is not a  Fredholm operator}\};$$
We denote by $\sigma_{lf}(T)$ and $\sigma_{uf}(T)$,   respectively the lower and upper semi-Fredholm spectra which are defined in the
same manner.

\begin{thh}\label{aze}
Let $(T(t))_{t\geq 0}$ be a $C_0$ semigroup on $X$ with infinitesimal generator $A$. Suppose that there exists $t_0>0$ such that $AT(t_0)$ is bounded. Then

$$\sigma_{uf}(T(t))\backslash \{0\}=\{e^{\lambda s}, \lambda\in\sigma_{uf}(A)\}$$
If $X$ is reflexive, we have :
$$\sigma_{lf}(T(t))\backslash \{0\}=\{e^{\lambda s}, \lambda\in\sigma_{lf}(A)\}$$
$$\sigma_{e}(T(t))\backslash \{0\}=\{e^{\lambda s}, \lambda\in\sigma_{e}(A)\}$$
\end{thh}

\p
Let $\mu=e^{\lambda t}\notin\sigma_{uf}(T(t))$, then $\dim N(e^{\lambda t}-T(t))<\infty$ and $R(e^{\lambda t}-T(t))$ is  closed. According \cite[Lemma 2.3]{taj}, $R(\lambda -A)$ is  closed. Since $N(\lambda-A)\subseteq N(e^{\lambda t}-T(t))$, we have $\dim N(\lambda-A)<\infty$.  Then $\lambda -A$ is upper semi-Fredholm. Therefore
$$\{\displaystyle e^{\lambda t}, \lambda\in\sigma_{uf}(A)\}\subseteq \sigma_{uf}(T(t))\backslash \{0\}\quad\quad (1).$$
Now, let $\mu\in\sigma_{uf}(T(t))\backslash\{0\}$ and $F=N(\mu-T(t))$, then $F$ is an  $A$-invariant and $T(t)$-invariant closed  subspace of $X$.
From lemma \ref{llla},  $A_{|F}$ is a bounded operator that generates the  $C_0$ semigroup $(T(t)_{|F})_{t\geq 0}$. By $(1)$, we have
$$\{\displaystyle e^{\lambda t}, \lambda\in\sigma_{uf}(A_{|F})\} \subseteq \sigma_{uf}(T(t)_{|F})\backslash \{0\} \quad\quad.$$
Since $A_{|F}$ is bounded, then $\sigma_{uf}(A_{|F})\neq \emptyset$. Let $\lambda_0\in\sigma_{uf}(A_{|F})$, then $\displaystyle e^{\lambda_0 t}\in \sigma_{uf}(T(t)_{|F})\subseteq \sigma(T(t)_{|F})=\{\mu\}$. Then $\mu=\displaystyle e^{\lambda_0 t} \in \{ e^{\lambda t}, \lambda\in\sigma_{uf}(A_{|F})\}$.\\
Now, let's show that $\lambda_0\in\sigma_{uf}(A)$. If $\lambda_0\notin\sigma_{uf}(A)$, so $A-\lambda_0$ is upper semi Fredholm, then
$\dim N(\lambda_0-A)<\infty$ and $R(\lambda_0-A)$ is closed. Since $N(\lambda_0-A_{|F})\subseteq N(\lambda_0-A)$, then
$\dim N(\lambda_0-A_{|F})<\infty$. Let $(\lambda_0-A_{|F})x_n\longrightarrow y$ with $(x_n)_n\subseteq F$. Since $R(\lambda_0-A)$ is closed, then  there exists  $x\in D(A)$ such that
$y=(\lambda_0-A)x$.
$B_{\lambda_0}(t)$ is bounded, so that then
$B_{\lambda_0}(t)(\lambda_0-A_{|F})x_n\longrightarrow  B_{\lambda_0}(t)(\lambda_0-A)x$, then $(\mu-T(t))x_n\longrightarrow (\mu-T(t))x$. Hence,
$(\mu-T(t))x=0$, consequently $x\in F$ and $y \in R(\lambda_0-A_{|F})$, which is absurd. Thus
$$\sigma_{uf}(T(t))\backslash \{0\}=\{ e^{\lambda t}, \lambda\in\sigma_{uf}(A)\}.$$
We now turn to the second statement.
Since $A$ is the generator of the  $C_0$ semigroup $(T(t))_{t\geq 0}$, then $A^*$ is the generator of the $C_0$ semigroup $(T^*(t))_{t\geq 0}$. Then
$$\sigma_{uf}(T^*(t))\backslash \{0\}=\{ e^{\lambda t}, \lambda\in\sigma_{uf}(A^*)\}.$$

By duality, we have  $$\sigma_{lf}(T(t))\backslash \{0\}=\left\{ e^{\lambda t}, \lambda\in\sigma_{lf}(A)\right\}$$
and consequently
$$\sigma_{e}(T(t))\backslash \{0\}= \{e^{\lambda t}, \lambda\in\sigma_{e}(A)\}. \blacksquare$$


Recall that  the ascent of $A$ is defined as $\mbox{asc}(A)=\mbox{inf}\{n\in \mathbb{N}:N(A^n)=N(A^{n+1})\}$
and the descent of $A$ is defined as $\mbox{des}(A)=\mbox{inf}\{n\in \mathbb{N}: R(A^n)=R(A^{n+1})\}$. Where $\mbox{inf } \emptyset=\infty$.
 We say that a closed operator $A$ is upper semi-Browder if it is
upper semi-Fredholm and has finite ascent.
Similarly, $A$ is lower semi-Browder if it is lower semi-Fredholm and has finite
descent.  $A$ is Browder if it is both lower and upper semi-Browder.
Equivalently, this means that $A$ is Fredholm and has finite both ascent and descent.

The Browder spectrum of $T$  is defined  by :
$$\sigma_{b}(T)=\{\lambda\in \mathbb{C}:  T-\lambda I \mbox{  is not a  Browder operator}\}.$$
Also, we denote by $\sigma_{lb}(T)$ and $\sigma_{ub}(T)$,   respectively the lower and upper semi-Browder spectra  which are defined in the
same manner. 

The following corollary is an immediate consequence of the previous theorem.

\begin{co}\label{sd}
Let $(T(t))_{t\geq 0}$ be a $C_0$ semigroup on  $X$ with infinitesimal generator $A$. Suppose that there exists $t_0>0$ such that $AT(t_0)$ is bounded and $X$ is reflexive.  Then

$$\sigma_{ub}(T(t))\backslash \{0\}=\{e^{\lambda t}, \lambda\in\sigma_{ub}(A)\}$$
If $X$ is reflexive, we have :
$$\sigma_{lb}(T(t))\backslash \{0\}=\{e^{\lambda t}, \lambda\in\sigma_{lb}(A)\}$$
$$\sigma_{b}(T(t))\backslash \{0\}=\{e^{\lambda t}, \lambda\in\sigma_{b}(A)\}$$

Since every differentiable $C_0$ semigroup matches the condition of the Theorem \ref{aze}, we have the following result.
\end{co}
\begin{co}\label{2}
Let $(T(t))_{t\geq 0}$ ba a differentiable $C_0$ semigroup on $X$ with infinitesimal generator $A$. Then
$$\sigma_{\ast}(T(t))\backslash \{0\}=\{e^{\lambda s}, \lambda\in\sigma_{\ast}(A)\}$$
with $\sigma_{\ast} \in \{\sigma_{uf}, \sigma_{lf}, \sigma_{f}, \sigma_{ub}, \sigma_{lf}, \sigma_{b} \}.$

\end{co}

\begin{remark}\label{ass}
Note that the  inclusion $\{\displaystyle e^{\lambda t}, \lambda\in\sigma_{*}(A)\}\subseteq \sigma_{*}(T(t))\backslash \{0\}$,  where $\sigma_{*}=\sigma_{e}, \sigma_{uf}, \sigma_{lf}, \sigma_{ub}, \sigma_{lb}, \sigma_{b}$ holds in the general setting where $T(t)$ is a $C_0$ semigroup.
This  inclusion  is strict as shown in the following example.
\end{remark}

\begin{ex}
Consider the translation group on the space $C_{2\pi}(\mathbb{R})$ of all $2\pi$
periodic continuous functions on $\mathbb{R}$ and denote its generator by $A$ (see \rm\cite[Paragraph I.4.15]{engel}). From \rm\cite[Examples 2.6.iv]{engel} we have, $\sigma(A) = i\mathbb{Z}$, then $e^{t\sigma(A)}$ is at most countable, therefore $e^{t\sigma_{*}(A)}$ are also. The spectra of the operators $T(t)$ are always contained in $\Gamma=\{z\in \mathbb{C} : |z|=1\}$ and contain the eigenvalues $e^{ikt}$ for $k \in \mathbb{Z}$. Since $\sigma (T(t))$ is closed, it follows from \rm\cite[Theorem IV.3.16]{engel}  that $\sigma (T(t)) = \Gamma$ whenever $t/2\pi \notin \mathbb{Q}$, then $\sigma (T(t))$ is not countable, so $\sigma_{*} (T(t))$ are also. Therefore the inclusions of the Remark \ref{ass} are strict. According Corollary \ref{2}, the semigroup $(T(t))_{t\geq 0}$ is not differentiable.
\end{ex}

\begin{ex}
On $C([0,1])$ we consider the operator $B$ given by
$$Bu:=mu''+qu'$$
where $m(0)=m(1)=0$, $\sqrt{m}\in C^1[0,1]$, $q\in C[0,1]$, and $q/ \sqrt{m}$ is bounded in $]0,1[$.
The operator $(B, D(B))$ generates an analytic  semigroup $(T(t))_{t\geq 0}$ of angle $\frac{\pi}{2}$ in $C[0,1]$, see \rm\cite[Theorem 4.20]{Pazy}, where $D(B):=\{u\in C[0,1]\cap C^2(0,1): Bu\in C[0,1]\}$. Since $(T(t))_{t\geq 0}$ is an analytic semigroup, $(T(t))_{t\geq 0}$ matches the condition of the Theorem \ref{aze}, then
$$\sigma_{b}(T(t))\backslash \{0\}=\{e^{\lambda s}, \lambda\in\sigma_{b}(B)\};$$
$$\sigma_{e}(T(t))\backslash \{0\}=\{e^{\lambda s}, \lambda\in\sigma_{e}(B)\}.$$
\end{ex}

Next, let $A$ be a closed  linear operator, $A$ is said to have the single
valued extension property at $\lambda_{0}\in\mathbb{C}$ (SVEP) if
for every  open neighborhood   $U\subseteq \mathbb{C}$ of
$\lambda_{0}$, the only  analytic function  $f: U\longrightarrow D(A)$ which satisfies
 the equation $(A-zI)f(z)=0$ for all $z\in U$ is the function $f\equiv 0$. $A$ is said to have the SVEP if $A$ has the SVEP for
 every $\lambda\in\mathbb{C}$.

Recall that a  closed operator $A$ is said to be
semi-regular if $R(A)$ is closed
and $N(A)\subseteq R^{\infty}(A)$,  where  $R^{\infty}(A)=\bigcap_{n\geq0}R(A^n)$. The semi-regular spectrum  of $A$ is defined by:
$$\sigma_{k}(A)=\{\lambda\in\mathbb{C}: A-\lambda I\mbox{  is not semi-regular }\}.$$

We have the following result.

\begin{pro}\label{gf}
Let $(T(t))_{t\geq 0}$ be a $C_0$ semigroup on $X$ with infinitesimal generator $A$. Then

$$\{e^{\lambda s}, \lambda\in\sigma_{k}(A)\} \subseteq\sigma_{k}(T(t))\backslash \{0\}\subseteq\{e^{\lambda s}, \lambda\in\sigma_{ap}(A)\}$$
Furthermore, if $A^*$ has the SVEP, we have
$$\{e^{\lambda s}, \lambda\in\sigma_{k}(A)\}=\sigma_{k}(T(t))\backslash \{0\}$$
\end{pro}

\p

According to \cite[Theorem 2.2]{Mulk}, we have the inclusion $\{e^{\lambda s}, \lambda\in\sigma_{k}(A)\}\subseteq\sigma_{k}(T(t))\backslash \{0\}$.\\
Now, let $\mu\in\sigma_{k}(T(t))\backslash\{0\}$ and $F=N(\mu-T(t))$, then $F$ is an $A$-invariant and $T(t)$-invariant closed  subspace of $X$.
From lemma \ref{llla},  $A_{|F}$ is a bounded operator  that generates the  $C_0$ semigroup $(T(t)_{|F})_{t\geq 0}$.  We have
$$\{\displaystyle e^{\lambda t}, \lambda\in\sigma_{k}(A_{|F})\}\subseteq \sigma_{k}(T(t)_{|F})\backslash \{0\} \quad\quad.$$
Since $A_{|F}$ is bounded, then $\sigma_{k}(A_{|F})\neq \emptyset$. Let $\lambda_0\in\sigma_{k}(A_{|F})$, then $\displaystyle e^{\lambda_0 t}\in \sigma_{k}(T(t)_{|F})=
\sigma(T(t)_{|F})=\{\mu\}$. Then $\mu=\displaystyle e^{\lambda_0 t} \in \{ e^{\lambda t}, \lambda\in\sigma_{k}(A_{|F})\}\cup\{0\}$.
Now, let us show that $\lambda_0\in\sigma_{*}(A)$. If $\lambda_0\notin\sigma_{*}(A)$, so $\lambda_0\notin\sigma_{*}(A_{|F})$, then $\lambda_0\notin\sigma_{k}(A_{|F})$, absurd. $\blacksquare$

The following corollary is an immediate consequence of the previous proposition.
\begin{co}
Let $(T(t))_{t\geq 0}$ be a $C_0$ semigroup on $X$ with infinitesimal generator $A$. 
If $X$ is reflexive and  $A$ has the SVEP, then
$$\{e^{\lambda s}, \lambda\in\sigma_{k}(A)\}=\sigma_{k}(T(t))\backslash \{0\}$$
\end{co}

\section{Differentiable $C_0$ semigroup}

\begin{lemma}\label{er}

Let $(T(t))_{t\geq 0}$ ba a differentiable $C_0$ semigroup on $X$ with infinitesimal generator $A$. Let $F=N(\mu-AT(t))$, $\mu\in\sigma(AT(t))\backslash\{0\}$. Then $(T(t)_{|F})_{t\geq 0}$ is   a uniformly continuous $C_0$-semigroup on $X$ with infinitesimal generator $A_{|F}$.

\end{lemma}

\p
For all $x\in F$, \begin{eqnarray*}
                    (\mu -AT(t))T(s)x &=& \mu T(s)x-AT(t)T(s)x \\
                                      &=& T(s)(\mu-AT(t))x=0
                  \end{eqnarray*}

For all $x\in F\cap D(A)$, \begin{eqnarray*}
                    (\mu -AT(t))Ax &=& \mu Ax-AT(s)Ax \\
                                      &=& A(\mu-AT(t))x=0
                  \end{eqnarray*}
Hence $F$ is a closed subspace of $X$,  $A$-invariant  and $T(t)$-invariant.
$\overline{D(A_{\mid F})}=F$. Indeed: Let $x\in F$ and $x_p=p\displaystyle\int^{t+\frac{1}{p}}_t T(t)xdx$ $p\geq 1$, then $x_p\in D(A)\cap F=D(A_{|F})$ and $x_p\longrightarrow x$,
 $p\rightarrow +\infty$. Then $x\in \overline{D(A_{|F})}$, so $\overline{D(A_{\mid F})}=F$.
Let us show that $A_{|F}$ is bounded. Let $\lambda_0\in\rho(A_{|F})$, then $R(\lambda_0,A_{|F})$ is bounded below, indeed: If not, there exists $(x_p)_p\subset F$  such that $\parallel  x_p \parallel=1$ and $R(\lambda_0,A_{|F})x_p\rightarrow 0$,  where $p\rightarrow +\infty$. Since  $(T(t))_{t\geq 0})$ is differentiable,  then  $AT(t_0)$ is  bounded, it is the same for $(\lambda_0-A)AT(t)$.
$(\lambda_0-A)AT(t)R(\lambda_0, A_{\mid F})x_p=AT(t)x_p\rightarrow 0$,  where $p\rightarrow +\infty$. So, $T(t_0)x_n\rightarrow 0$,  where $n\rightarrow +\infty$, then  $x_p\rightarrow 0$ which contradicts the fact that $\parallel  x_n \parallel=1$. Hence, there exists $C>0$ such that for all $x\in F$ $\parallel R(\lambda_0, A_{\mid F})x\parallel\geq C\parallel x\parallel$. Then, for all $x\in D(A_{\mid F})$
\begin{center}
    $\parallel A x\parallel\leq (|\lambda_0|+ C^{-1})\parallel x \parallel$
\end{center}
Hence the result. $\blacksquare$

\begin{thh}

Let $(T(t))_{t\geq 0}$ ba a differentiable $C_0$ semigroup on $X$ with infinitesimal generator $A$. Then for $t>0$, we have
$$\sigma(AT(t))\cup \{0\}=\{\lambda e^{\lambda t},\lambda\in\sigma(A)\cup\{0\}\}$$
$$\sigma_{ap}(AT(t))\cup \{0\}=\{\lambda e^{\lambda t},\lambda\in\sigma_{ap}(A)\cup\{0\}\}$$
\end{thh}
\p
From \cite[Lemma 4.6]{Pazy},  we have  $\{\lambda e^{\lambda t},\lambda\in\sigma(A)\}\subseteq \sigma(AT(t)).$
Let $\mu\in\sigma(AT(t))\backslash \{0\}$. From lemma \ref{er}
 $(T(t)_{|F})_{t\geq 0}$ is   a uniformly continuous $C_0$-semigroup on $X$ with infinitesimal generator $A_{|F}$ where $F=N(\mu-AT(t))$, $\mu\in\sigma(A(t))\backslash\{0\}$. Then, $T(t)_{|F}=e^{t A_{|F}}$. Since the function $f: z\longrightarrow ze^{tz}$ is analytic  on $\mathbb{C}$, from the spectral mapping theorem we have:
$$\sigma((AT(t))_{|F})=\sigma(A_{|F}T(t)_{|F})=\sigma(A_{|F} e^{t A_{|F}})=\sigma(f(A_{|F}))=f(\sigma(A_{|F})).$$
$A_{|F}$ is bounded, then $\sigma(A_{|F})\neq \emptyset$. Let $\lambda_0\in\sigma(A_{|F})$, then $\lambda_0e^{t\lambda_0}=f(\lambda_0)\in\sigma((AT(t))_{|F})=\{\mu\}$. Therefore, $\mu=\lambda_0e^{ t\lambda_0}$. Now, let us show that $\lambda_0\in\sigma(A)$.
 If $\lambda_0-A$ is bijective, so $\lambda_0-A_{|F}$ is injective. Let $y\in F$, $\lambda_0-A$ is surjective, then there exists $x\in D(A)$ such that
 $y=(\lambda_0-A)x$, since $y\in F$ then $(\lambda_0 e^{\lambda_0 t}-AT(t))y=(\lambda_0-A)B'_{\lambda_0}(t)x=0$, so $B'_{\lambda_0}(t)y=0$
 and as a result $(\lambda_0 e^{\lambda_0 t}-AT(t))y=(\lambda_0-A)B'_{\lambda_0}(t)x=0$, thus $x\in F$. Hence, $\lambda_0\in\rho(A_{|F})$, absurd.
 Therefore,  $$\sigma(AT(t))\cup \{0\}=\{\lambda e^{\lambda t},\lambda\in\sigma(A)\cup\{0\}\}.$$
 Now, let $\lambda_0\in \sigma_{ap}(A)$, then  there exists $(x_p)_p\subset D(A)$ such that $\parallel x_p \parallel=1$ and $(\lambda_0-A)x_p\longrightarrow 0$ where $p\rightarrow +\infty$. We have  $B'_{\lambda_0}(t)$ is bounded, then
  $(\lambda_0e^{\lambda_0 t}-AT(t))x_p=B'_{\lambda_0}(t)(\lambda_0-A)x_p\longrightarrow 0$ where $p\rightarrow +\infty$. Therefore $\lambda_0e^{\lambda_0 t} \in\sigma_{ap}(AT(t)).$  So, we have
$$\{\lambda e^{\lambda t},\lambda\in\sigma_{ap}(A)\cup\{0\}\} \subseteq\sigma_{ap}(AT(t))\cup \{0\}$$
Let $\mu\in\sigma_{ap}(AT(t))\backslash\{0\}$. Since $\sigma_{ap}((AT(t))_{|F})=\sigma(AT(t)_{|F})=f(\sigma_{ap}(A_{|F}) )=\{\mu\}$.
$A_{|F}$ is bounded, then $\sigma_{ap}(A_{|F})\neq \emptyset$. Let $\lambda_0\in \sigma_{ap}(A_{|F})$, then $\mu=\lambda_0e^{\lambda_0 t}$.
We have $\lambda_0\in\sigma_{ap}(A)$. Indeed: suppose that $\lambda_0\notin\sigma_{ap}(A)$, then $\lambda_0-A$ is injective  with closed range, so
$\lambda_0-A_{|F}$ is injective. Let $(\lambda_0-A_{|F})x_n\longrightarrow y$. Since $R(\lambda_0-A$ is closed then  there exists $x\in D(A)$ such that $y=(\lambda_0-A_{|F})x $. $B'_{\lambda_0}(t)$ is bounded, then $B'_{\lambda_0}(t)(\lambda_0-A_{|F})x_n\longrightarrow B'_{\lambda_0}(t)(\lambda_0-A_{|F})x=(\lambda_0e^{\lambda_0 t}-AT(t)x)x$. $x_n\in F$, then  $B'_{\lambda_0}(t)(\lambda_0-A_{|F})x_n=0$. Then $x\in F$. Therefore $\lambda_0\in\rho_{ap}(A_{|F})$, absurd. Consequently, $\sigma_{ap}(AT(t))\cup \{0\}=\{\lambda e^{\lambda t},\lambda\in\sigma_{ap}(A)\cup\{0\}\}$. $\blacksquare$

\begin{co}\label{aqs}

Let $(T(t))_{t\geq 0}$ ba a differentiable $C_0$ semigroup on $X$ with infinitesimal generator $A$. Then for $t>0$, we have
$$\forall n\in\mathbb{N}, \forall t>0\quad\quad\quad\sigma_*(T(t)^{(n)})\cup \{0\}=\{\lambda^n e^{\lambda t},\lambda\in\sigma_*(A)\cup\{0\}\}$$
where $\sigma_*=\sigma, \sigma_{ap}$
\end{co}

\p Let $n\in\mathbb{N}$, $t>0$. Since $(T(t))_{t\geq 0}$ is a differentiable $C_0$ semigroup, then $T(t)^{(n)}$ is a bounded operator, furthermore
$$T(t)^{(n)}= A^n T(t)=\left(AT\left(\frac{t}{n}\right)\right)^n.$$
From the spectral mapping theorem, we have
\begin{eqnarray*}
  \sigma_*(T(t)^{(n)})\cup \{0\} &=&  \sigma_*\left(AT\left(\frac{t}{n}\right)\right)^n\cup \{0\}=\left\{\mu^n: \mu\in\sigma_*\left(AT\left(\frac{t}{n}\right)\right)\right\}\{0\} \\
   &=& \left\{\left(\lambda e^{\lambda \frac{t}{n}}\right)^n,\lambda\in\sigma_*(A)\cup\{0\}\right\} \\
   &=& \{\lambda^n e^{\lambda t},\lambda\in\sigma_*(A)\cup\{0\}\}.\blacksquare
\end{eqnarray*}

\begin{ex}
Let $X := \ell^p$ and $M_q$ be the multiplication operator defined by  $M_q(x_n)_{n\in \mathbb{N}}=(q_nx_n)_{n\in \mathbb{N}}$ with $q=(q_n)_{n\in \mathbb{N}}$ and $q_n=-n+i n^2$. Then $\sigma(M_q) = \{q_n : n \in \mathbb{N}\}$. From \rm\cite[Counterexamples II.4.16]{engel} the semigroup $(T_q(t))_{t\geq 0}$ generated by $M_q$ is differentiable.  By Corollary \ref{aqs}, for $t> 0$, we have $\sigma(T_q(t)^{(p)})\cup \{0\}=\{(-n+i n^2)^pe^{t (-n+i n^2)} : n\in \mathbb{N}\}$.
\end{ex}

\begin{ex}
Consider the heat equation in $L^p(0, \pi)$.\\

$\left\{
  \begin{array}{ll}
    \displaystyle\frac{\partial u}{\partial t}(t,x)= \displaystyle\frac{\partial^2 u}{\partial x^2}(t,x),(t,x)\in\mathbb{R}^+ \times (0,\pi) & \hbox{} \\
    u(t,0)=0=u(t,\pi), t\geq 0 & \hbox{} \\
    u(0,x)=f(x)\quad x\in(0,\pi) & \hbox{.}
  \end{array}
\right.$

Let $p>  2$. On $X = L^p(0, \pi)$ consider
the operator defined by
$$Af(x)=f''(x)$$
 with domain $D(A)=W^{2,p}(0,\pi)\cap W_0^{1,p}(0,\pi)$, $x\in(0,\pi)$
 where $$W_0^{1,p}=\{f\in W^{1,p}(0,\pi): f(0)=0=f(\pi)\}.$$

 The operator $A$  is self-adjoint. The corresponding  semigroup $(T(t))_{t\geq 0}$ is analytic, so is differentiable.

For each $f\in W^{2,p}(0,\pi)\cap W_0^{1,p}(0,\pi)$ the unique solution of the equation is given by $$u(t, x) = (T(t)f)(x).$$

The spectrum of $A$ is $\sigma(A)=\{-n^2; n\geq 1\}.$  Since $int(\sigma(A))=\emptyset$, then $A$  and $A^*$ have SVEP. So, $\sigma(A)=\sigma_k(A)=\sigma_{ap}(A)=\sigma_{su}(A)$.

From Proposition \ref{gf}, we have $\sigma_{k}(T(t))\cup\{0\}=\{ e^{-tn^2}, n\geq 1\}\cup\{0\}$.

Also, by Corollary \ref{aqs}, we have  $$\sigma(T(t)^{(p)})\cup\{0\}=\sigma_{ap}(T(t)^{(p)})\cup\{0\}=\{ (-1)^pn^{2p}e^{-tn^2}, n\geq 0\}.$$
\end{ex}


\begin{thebibliography}{99}




\bibitem{Mulk} \textsc{A. Elkoutri and M. A. Taoudi,} \emph{Spectre singulier pour les g\'{e}n\'{e}rateurs
des semi-groupes }, C. R. Acad. Sci. Paris Sér. I Math., 333 (2001), no. 7, 641-644
(French).

\bibitem{engel}  \textsc{K. J. Engel and R. Nagel,} \emph{One-Parameter Semigroups for Linear Evolution Equations},
Graduate Texts in Mathematics, vol. 194, Springer-Verlag, New York, 2000.




\bibitem{Pazy}    \textsc{A. Pazy,} \emph{Semigroups of Linear Operators and Applications to Partial Differential
Equations}, Applied Mathematical Sciences, vol. 44, Springer-Verlag,  New York 1983.

\bibitem{taj1} \textsc{A. Tajmouati and H. BOUA}  \emph{ Spectral Mapping Theorem for $C_0$-Semigroups of Drazin Spectrum, } Bol. Soc. Paran. Mat., to appear.


\bibitem{taj} \textsc{A. Tajmouati, M. Amouch,  M.R.F. Alhomidi Zakariya}  \emph{Spectral Inclusion for $C_0$-Semigroups Drazin Invertible and Quasi-Fredholm Operators, }
Bull. Malays. Math. Sci. Soc., to appear.

\bibitem{taj2} \textsc{A. Tajmouati, M. Amouch,  M.R.F. Alhomidi Zakariya}  \emph{ Spectral Equality for $C_0$-Semigroups and
Spectral Inclusion of B-Fredholm, } Rend. Circ. Mat. Palermo,  65(3) (2016), 425-434.


\end{thebibliography}
\end{document}